\documentclass[aip,amsmath,preprint,reprint]{revtex4-1}
\usepackage{hyperref}
\usepackage{graphicx}

\newcommand{\de}{\,\mathrm{d}}

\begin{document}

\title{Global approximation to the Boys functions for vectorized computation}

\author{Dimitri N. Laikov}
\email[]{laikov@rad.chem.msu.ru}
\affiliation{Chemistry Department, Moscow State University, 119991 Moscow, Russia}

\date{\today}

\begin{abstract}
A fast approximation to the Boys functions
(related to the lower incomplete gamma function of half-integer parameter)
by a single closed-form analytical expression for all argument values
have been developed and tested.
Besides the exponential function needed anyway for downward recursion,
it uses a small number of addition, multiplication, division, and square root operations,
and thus is straightforward to vectorize.
\end{abstract}

\maketitle

\section{Introduction}

The Boys~\cite{B50} functions
\begin{equation}
\label{eq:fn}
F_n(x) = \int\limits_0^1 t^{2n} \exp\left(-xt^2\right) \de t,
\end{equation}
related to the lower incomplete gamma function of half-integer parameter,
are the only special functions in the analytical computation
of one- and two-electron integrals~\cite{MD78} over Gaussian-type basis functions
in molecular electronic structure theory ---
they have to be evaluated many times, so that fast and accurate approximations are needed.
The most traditional approach is to use a piecewise polynomial approximation~\cite{B50,MD78,GJP91,WO15,PK20} ---
this may well be the fastest serial method, but table lookup does not seem to be vectorizable,
and vectorization seems now to be the only way forward toward faster computation.
A global analytic approximation to functions~(\ref{eq:fn}) would be an ideal solution,
historically there is an early work~\cite{SM71} on rational approximation
followed by another dramatic attempt~\cite{FO94} more than 20 years later,
and 30 years of silence since then --- should this mean any further work is hopeless?
We have found a better functional form of the approximation,
first for the error function~\cite{L25a}, and hence $F_0(x)$ as well,
and later generalized it to $n>0$ and report it here.

\section{Methodology}

\begin{figure}[h]
\includegraphics[scale=1.1]{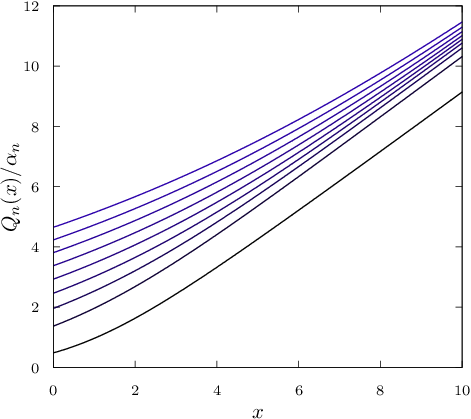}
\caption{\label{fig:1} Functions $Q_n(x)$, $n=0,\dots,8$ (black to blue).}
\end{figure}

Our new functional form
\begin{equation}
\label{eq:fq}
F_n(x) = \frac{c_n}{\left(x + \bigl( Q_n(x) \bigr)^{n+1/2} \exp(-x) \right)^{n+1/2}}
\end{equation}
defines $F_n(x)$ in terms of the new functions $Q_n(x)$,
\begin{equation}
c_n = \frac{(2n)!\; \sqrt{\pi}}{n!\; 2^{2n+1}},
\end{equation}
and the long-range limit
\begin{equation}
\label{eq:limf1}
\lim_{x\to\infty} F_n(x) = \frac{c_n}{x^{n+1/2}}
\end{equation}
is reached with exponential convergence, as it should.
The functions have the linear asymptote
\begin{equation}
\label{eq:qa}
\lim_{x\to\infty} Q_n(x) = \alpha_n x,
\end{equation}
\begin{equation}
\label{eq:qb}
\lim_{x\to\infty} Q_n(x) - \alpha_n x = \beta_n + O\left(x^{-1} \right),
\end{equation}
\begin{equation}
\label{eq:aq}
\alpha_n = \bigl(2 c_{n+1} \bigr)^{-2/(2n+1)},
\end{equation}
\begin{equation}
\label{eq:bq}
\beta_n = \frac{2n-1}{2n+1} \alpha_n,
\end{equation}
while at $x=0$ the value and the first derivative are
\begin{equation}
\label{eq:q0}
Q_n(0) = \bigl((2n+1) c_n \bigr)^{4/(2n+1)^2},
\end{equation}
\begin{equation}
\label{eq:q1}
Q'_n(0) = \frac{2 Q_n(0)}{2n+1}\left(\frac{2n+5}{2n+3} - \bigl((2n+1) c_n \bigr)^{-2/(2n+1)} \right),
\end{equation}
and higher derivatives have ever lengthier expressions.
Fig.~\ref{fig:1} shows the functions to be remarkably regular!

The rational approximation
\begin{equation}
\label{eq:qr}
\tilde{Q}_n(x) = \frac{\sum_{k=0}^{N+1} A_{kn} x^k}{1 + \sum_{k=1}^{N} B_{kn} x^k} \approx Q_n(x)
\end{equation}
is a very natural choice as it can match both the long- and short-range behavior
(\ref{eq:qa}), (\ref{eq:qb}), (\ref{eq:q0}), (\ref{eq:q1}),
and even the higher terms of the series both at $x=0$ and $x\to\infty$.

Numerical evaluation of the rational function~(\ref{eq:qr}) in finite-precision arithmetic
may suffer (see below) from rounding errors,
but it can be done in mathematically equivalent form
\begin{equation}
\label{eq:pqz}
\tilde{Q}_n(x) = q_{0n} + q_{1n} x + q_{2n} \frac{u_n(x)}{v_n(x)} x^2,
\end{equation}
$q_{0n} \equiv Q_n(0)$ as in~(\ref{eq:q0}) and $q_{1n} \equiv Q'_n(0)$ as in~(\ref{eq:q1}),
where the polynomials
\begin{equation}
\label{eq:pz}
u_n(x) = \prod_{k = 1}^{L} \left(\left(x - y_{kn} \right)^2 + a_{kn}\right)
 \prod_{k = L + 1}^{N - L - 1} \left(x - y_{kn} \right),
\end{equation}
\begin{equation}
\label{eq:qz}
v_n(x) = \prod_{k = 1}^{M} \left(\left(x - z_{kn} \right)^2 + b_{kn}\right)
 \prod_{k = M + 1}^{N - M} \left(x - z_{kn} \right),
\end{equation}
are evaluated as written in terms of their real and complex roots,
with much less loss of accuracy and the same number of operations.
The product of monomials and binomials with negative-real-part roots
can be further gathered into one polynomial factor.

We prefer the approximations with uniform relative accuracy,
minimizing the error measure
\begin{equation}
\label{eq:e}
E_n = \max\limits_{x \ge 0} \left|\frac{\tilde{F}_n(x)}{F_n(x)} - 1 \right|
\end{equation}
with respect to the parameters $\{A_{kn}\}$ and $\{B_{kn}\}$,
with the four constraints (\ref{eq:qa}), (\ref{eq:qb}), (\ref{eq:q0}), (\ref{eq:q1}) built into them.

Downward recursion
\begin{equation}
\label{eq:r}
F_{m-1}(x) = \frac{2x F_m(x) + \exp(-x)}{2m + 1}
\end{equation}
is known to be stable and yields the whole set of values for $m=0,\dots,n$
as needed in practice, the same value of $\exp(-x)$ is also used in our approximation~\ref{eq:fq}.

The helpful cut-off values~$z_{nb}$ as solutions of
\begin{equation}
\frac{z_{nb}^{n+1/2}}{c_n} F_n\bigl(z_{nb}\bigr) = 1 - 2^b
\end{equation}
can be tabulated for $b$-bit precision 
and used to switch to the asymptotic formula~\ref{eq:limf1}
for $x > z_{nb}$.

\section{Computations}

We have written a C code~\cite{L25a} to optimize the uniform approximations using extended precision arithmetic,
doing most work in 256-bit precision, and then running it once again,
for the love of the art, in 512-bit precision.

For all $n=0,\dots,36$ we have sought solutions from $N=2$ onwards until reaching
the accuracy of 64 bits, $-\log_2 E_n \approx 64$, leading to $N\le 20$.
In some few $(n,N)$ cases we could only find solutions having $B_{Nn} < 0$
and had to reject them.

We have gathered 547 solutions and formatted them as a C file
(see \href{ftp://}{supplementary material}) and added more code
to get a program for thorough testing of the approximations
working in 24-bit (mantissa) single, 53-bit double, 64-bit long double, and 113-bit quadruple precision.
At first, we were scared to discover an unexpectedly high loss of accuracy on going to higher $n$,
understanding its cause to be the positive values of polynomial coefficients $A_{kn}$ and $B_{kn}$ in (\ref{eq:qr}),
also meaning there must be polynomial roots with positive real parts.
With little belief, we did try Eqs.~(\ref{eq:pz}), (\ref{eq:qz}), and then~(\ref{eq:pqz}),
and our hopes were rewarded --- the errors have dropped to a low enough level.

\begingroup
\begin{table}
\caption{\label{tab:e}Approximation accuracy and rounding errors.}
\begin{ruledtabular}
\begin{tabular}{r|rrrrr|rrrrr}
 & & \multicolumn{3}{c}{53-bit double} &
 & & \multicolumn{3}{c}{24-bit single} & \\
$n$ & $N$ & $\bar{\epsilon}_0$ & $\bar{\epsilon}_{n-1}$ & $\bar{\epsilon}_n$ & $\epsilon_n$
    & $N$ & $\bar{\epsilon}_0$ & $\bar{\epsilon}_{n-1}$ & $\bar{\epsilon}_n$ & $\epsilon_n$ \\
\hline
  0 & 10 & 52.0 &      & 52.0 & 58.0 & 4 & 22.9 &      & 22.9 & 27.3 \\
  1 &    & 50.2 &      & 49.6 & 52.7 &   & 21.4 &      & 21.1 & 23.7 \\
  2 & 11 & 49.9 & 50.0 & 49.0 & 53.8 & 5 & 21.0 & 21.1 & 19.8 & 27.3 \\
  3 &    & 49.0 & 49.1 & 48.4 & 51.5 &   & 20.4 & 20.4 & 19.0 & 25.5 \\
  4 &    & 48.7 & 48.8 & 47.7 & 51.1 &   & 20.1 & 20.2 & 18.5 & 24.2 \\
  5 &    & 48.7 & 48.9 & 46.9 & 50.7 &   & 19.8 & 20.0 & 18.3 & 23.3 \\
  6 &    & 48.3 & 48.3 & 47.0 & 51.1 &   & 19.4 & 19.5 & 17.7 & 22.5 \\
  7 & 12 & 48.1 & 48.3 & 46.6 & 51.3 &   & 19.2 & 19.4 & 17.5 & 21.9 \\
  8 &    & 48.2 & 48.2 & 46.2 & 51.0 &   & 19.0 & 19.0 & 17.0 & 21.3 \\
  9 &    & 47.8 & 47.9 & 45.5 & 50.1 &   & 18.9 & 18.8 & 16.7 & 20.9 \\
 10 &    & 47.6 & 47.5 & 45.4 & 49.3 & 6 & 18.9 & 18.7 & 16.4 & 24.9 \\
 11 & 13 & 47.8 & 47.4 & 45.0 & 51.0 &   & 19.2 & 18.6 & 16.1 & 24.5 \\
 12 &    & 47.6 & 47.4 & 44.6 & 50.7 &   & 18.6 & 18.4 & 16.1 & 24.2 \\
 13 &    & 47.4 & 47.2 & 44.8 & 50.3 &   & 18.5 & 18.3 & 15.9 & 23.8 \\
 14 &    & 47.3 & 47.0 & 44.5 & 49.9 &   & 18.4 & 18.1 & 15.7 & 23.5 \\
 15 &    & 47.1 & 46.5 & 43.9 & 49.6 &   & 18.5 & 17.8 & 15.1 & 23.3 \\
 16 &    & 47.0 & 46.6 & 44.1 & 49.2 &   & 18.4 & 17.7 & 15.3 & 23.0 \\
 17 &    & 47.0 & 46.4 & 43.8 & 48.8 &   & 18.3 & 17.6 & 14.8 & 22.8 \\
 18 & 14 & 47.4 & 46.2 & 43.4 & 51.6 &   & 18.4 & 17.3 & 14.6 & 22.6 \\
 19 &    & 47.2 & 46.1 & 43.2 & 51.3 &   & 18.3 & 17.3 & 14.6 & 22.4 \\
 20 &    & 47.1 & 46.0 & 43.3 & 51.1 &   & 18.0 & 17.2 & 14.6 & 22.2 \\
 21 &    & 47.1 & 45.7 & 43.1 & 50.8 &   & 18.3 & 16.9 & 13.9 & 22.0 \\
 22 &    & 47.0 & 45.8 & 43.2 & 50.5 &   & 18.1 & 17.1 & 14.2 & 21.8 \\
 23 &    & 46.8 & 45.6 & 42.8 & 50.3 &   & 18.0 & 16.7 & 13.7 & 21.7 \\
 24 &    & 46.7 & 45.6 & 42.9 & 50.1 &   & 17.7 & 16.7 & 14.2 & 21.5 \\
 25 &    & 46.5 & 45.6 & 43.0 & 49.9 &   & 17.7 & 16.5 & 13.9 & 21.4 \\
 26 &    & 46.7 & 45.0 & 42.4 & 49.7 &   & 17.8 & 16.2 & 13.3 & 21.2 \\
 27 &    & 46.6 & 45.5 & 42.7 & 49.5 &   & 17.5 & 16.4 & 13.8 & 21.1 \\
 28 &    & 46.4 & 45.1 & 42.3 & 49.3 &   & 17.4 & 15.9 & 13.1 & 21.0 \\
 29 &    & 46.4 & 45.1 & 42.3 & 49.1 &   & 17.7 & 16.0 & 13.5 & 20.9 \\
 30 &    & 46.4 & 45.1 & 42.2 & 48.9 &   & 17.6 & 16.2 & 13.4 & 20.8 \\
 31 &    & 46.3 & 44.5 & 41.8 & 48.8 &   & 17.4 & 16.1 & 13.4 & 20.7 \\
 32 &    & 46.3 & 44.9 & 42.2 & 48.6 &   & 17.6 & 16.0 & 13.0 & 20.6 \\
 33 &    & 46.3 & 44.4 & 41.8 & 48.4 &   & 17.4 & 15.7 & 12.9 & 20.5 \\
 34 &    & 46.3 & 44.6 & 41.8 & 48.3 &   & 17.5 & 15.6 & 12.8 & 20.4 \\
 35 &    & 46.1 & 44.3 & 41.5 & 48.2 &   & 17.5 & 15.6 & 12.7 & 20.3 \\
 36 &    & 46.0 & 44.3 & 41.8 & 48.0 &   & 17.3 & 15.7 & 13.0 & 20.2 \\
\end{tabular}
\begin{flushleft}
The values of $E$ (exact arithmetics) and $\bar{E}$ (finite precision)
of Eq.~(\ref{eq:e}) are shown as negative binary logarithms: $\epsilon = -\log_2 E$.
\end{flushleft}
\end{ruledtabular}
\end{table}
\endgroup

Table~\ref{tab:e} shows the accuracy of approximations
with a choice of orders $N$ for working in either double or single precision,
such that the exact-arithmetic error is well buried under the rounding error.
For each $n=0,\dots,36$ we computed,
on a very dense ($\sim 2^{20}$) equally-spaced set of points $0\le x < z_{nb}$,
the approximate values of $F_n(x)$ and the full recursion~(\ref{eq:r}) to get $F_m(x)$, $m=0,\dots,n - 1$,
and the maximum relative errors in all values were determined.

The downward recursion~\ref{eq:r} is seen to recover more accurate $F_m(x)$ for $m<n$,
and this is good as higher accuracy is needed in practice for lower $m$.
We get $F_0(x)$ to no worse than 46 bits with double-
and 17 bits with single-precision arithmetic even for $n=36$.
We find this level of accuracy more than enough for most practical applications!

Further tricks can be used: to get $F_n(x)$ with a few more bits,
$F_{n+1}(x)$ can be computed instead followed by the downward recursion;
if the highest accuracy is needed for $m=0,\dots,\bar{n}$, $\bar{n}<n$,
two values $F_{\bar{n}}(x)$ and $F_n(x)$ can be computed by the best approximations,
and two recursions (the one from $\bar{n}$ down to $0$, and the other from $n$ down to $\bar{n}+1$) can be done.
Other choices can be made, for example,
full 24-bits of single-precision values can be better produced
by rounding the values computed in double precision with about 28-bit accurate approximations.

Our dataset of approximation coefficients can be used to automatically generate
the full code, with all loops unwound, for fast computation of the Boys functions~(\ref{eq:fn}),
that can also be vectorized in the same way as we did~\cite{L25a} before.

\section{Conclusions}

The new closed-form expression we have found for the global approximation of the Boys functions
works well --- a rather small number of polynomial terms is enough to reach a high accuracy,
finite-precision floating-point computations can be done with a low enough rounding error
and are well vectorizable.

\section{Data availability}

See \href{https://doi.org/insert.later/}{supplementary material} for all approximation coefficients
and a computer code to test them.

\bibliography{fn}

\clearpage

\end{document}